\documentclass[a4paper,english, oneside, reqno]{amsart}
\usepackage[a4paper]{geometry}
\usepackage{amsthm,amsmath,amssymb,mathtools}
\usepackage{enumerate,enumitem,csquotes,verbatim}
\usepackage{bm}
\usepackage{setspace}
\usepackage[numbers]{natbib}
\usepackage[colorlinks=true, linkcolor=blue, citecolor=blue]{hyperref}
\usepackage{cleveref}
\usepackage{color}
\usepackage{tikz}

\theoremstyle{plain}
\newtheorem*{theorem*}{Theorem}
\newtheorem{theorem}{Theorem}
\newtheorem{lemma}[theorem]{Lemma}

\newtheorem*{claim*}{Claim}

\theoremstyle{definition}

\theoremstyle{remark}

\def\N{\mathbb{N}}
\def\Z{\mathbb{Z}}

\title{Weak Dissolution of Groups for Rearrangement}

\author[P. Kamthorntaksina]{Pholphum Kamthorntaksina}
\address[Pholphum Kamthorntaksina]{
Department of Mathematics and Computer Science, Faculty of Science, Chulalongkorn University, Bangkok 10330, Thailand}
\email{pholphumk@gmail.com}
\author[T. Kittipassorn]{Teeradej Kittipassorn}
\address[Teeradej Kittipassorn]{
Department of Mathematics and Computer Science, Faculty of Science, Chulalongkorn University, Bangkok 10330, Thailand, and Centre of Excellence in Mathematics, Ministry of Higher Education, Science, Research and Innovation, Thailand}
\email{teeradej.k@chula.ac.th}
\author[P. Srikirinth]{Patid Srikirinth}
\address[Patid Srikirinth]{
Department of Mathematics and Computer Science, Faculty of Science, Chulalongkorn University, Bangkok 10330, Thailand}
\email{cirhale@gmail.com}
\author[P. Tulabadee]{Puttimet Tulabadee}
\address[Puttimet Tulabadee]{
Department of Mathematics and Computer Science, Faculty of Science, Chulalongkorn University, Bangkok 10330, Thailand}
\email{potter.potter2545@gmail.com}
\author[N. Witchayapongkul]{Nattaporn Witchayapongkul}
\address[Nattaporn Witchayapongkul]{
Department of Mathematics and Computer Science, Faculty of Science, Chulalongkorn University, Bangkok 10330, Thailand}
\email{natapon.tom2@gmail.com}

\allowdisplaybreaks

\begin{document}
\maketitle

\begin{abstract}
Suppose that $k$ of $n$ classrooms are to be dissolved and their students moved to
the other classrooms so that the numbers of students in any two remaining classrooms
differ by at most $d$.
We give a necessary and sufficient condition on the numbers of students for this to
be possible whichever $k$ classrooms are dissolved.
Moreover, we determine how small the largest classroom can be when this is possible.
In particular, when all classrooms have different sizes, the smallest possible size
of the largest classroom is $d+\lceil (n-d-1)(n-d)/(2k)\rceil$ if $k<(n-d-1)/2$, and
it is $n$ or $n+1$ otherwise.
This extends the results of Siriseriwan and Sripatak, who studied the case $d=0$.
\end{abstract}

\section{Introduction}

The almost equal partitioning, or balanced partitioning, problem is a classical topic that has been studied for many years.
One well-studied problem concerns graph partitioning, including edge partitioning and vertex partitioning into balanced components under various constraints, such as degree or coloring conditions; see, for example, \cite{Amir, Andreev, Bollobas3, Deligkas, Teeradej, Krauthgamer, Li}.

In this paper, we study a different type of balancing problem that is motivated by the following situation.
There are $n$ classrooms in a certain school.
 This semester, the school director would like to reduce the number of classrooms for better management by dissolving $k$ classrooms and redistributing the students from those classrooms among the remaining $n-k$ classrooms so that the numbers of students in the classrooms are as equal as
possible. The students in the remaining classrooms stay in their own classrooms, while
students from the same dissolved classroom may be moved to different classrooms. To be fair to the students, the $k$ classrooms to be dissolved will be chosen arbitrarily and the dissolution must be done for any choice of $k$ classrooms.

Ideally, the remaining classrooms would have exactly the same number of
students, but this is possible only if the total number of students is divisible
by $n-k$.
Therefore, we allow the numbers of students in any two remaining classrooms to
differ by at most $d$, where $d$ is a non-negative integer.
Note that such a rearrangement is always possible when $d$ is large enough, but
it may be impossible when $d$ is small.

Let us now formalize this situation.
Let $n\in\N$ be the number of classrooms, $k\in\{1, 2, \dots, n-1\}$ the number of classrooms to be dissolved, $d$ a non-negative integer, and $a_i$ the number of students in classroom $i$ for any $i\in\{1,2,\dots, n\}$.
For convenience, we always assume that $a_1 \leq a_2 \leq
\dots \leq a_n$.
We say that a collection of rooms is \emph{$d$-dissolvable} if their students
can be moved to the other rooms so that the numbers of students in any two
remaining rooms differ by at most $d$. A non-decreasing sequence $(a_1, a_2, \dots, a_n)\in\N^n$ of the numbers of students is \emph{$(k,d)$-dissolvable} if every $k$ rooms are $d$-dissolvable. When $d=0$, the remaining rooms must have exactly the same number of students.
Since we only require the numbers of students in the remaining rooms to differ
by at most $d$, we call such a rearrangement a \emph{weak dissolution}.

Our first aim is to characterize $(k,d)$-dissolvable sequences.

Since classrooms have limited capacity, the school cannot make its classrooms
arbitrarily large. On the other hand, if all classrooms have different sizes, then the smaller
remaining classrooms always need some extra students, and when every classroom is
small, the dissolved classrooms may not have enough students to provide them.
Hence, our second aim is to find, for fixed $n$, $k$ and $d$, the smallest
possible number of students in the largest classroom for which such a
rearrangement is possible.

Let $N(n,k,d)$ be the minimum value of $a_n$ over
all $(k,d)$-dissolvable sequences and $N'(n,k,d)$ the minimum value of $a_n$ over
all strictly increasing $(k,d)$-dissolvable sequences.
Note that $N(n,k,d)$ is at most $2$ by \Cref{thm:alpha}\ref{thm2.1}, and it
equals $1$ when $d\ge1$, since the sequence $(1,1,\dots,1)$ is $(k,d)
$-dissolvable.
Hence, the problem becomes nontrivial when the classrooms are required to have
different sizes.

In 2004, Siriseriwan~\cite{Wacharasak} characterized strictly increasing $
(k,0)$-dissolvable sequences and determined $N'(n,k,0)$ for $k=1,2$.
In the following year, Sripatak~\cite{Piyashat} extended both results to every
$1\le k\le n-1$ and showed that
\[
N'(n,k,0)=  \begin{cases}\left\lceil \frac{(n-1)n}{2k} \right\rceil & \text{if} \quad 1\le k< \frac{n-1}{2}, \\
    n & \text{if} \quad \frac{n-1}{2} \le k\le n-1 \text{ and } (n-k) \mid \frac{n(n+1)}{2}, \\
    n+1 & \text{if} \quad \frac{n-1}{2} \le k\le n-1 \text{ and } (n-k) \nmid \frac{n(n+1)}{2}.
\end{cases}
\]
In 2006, Thamrongthanyalak~\cite{Athipat} obtained the analogous results for
$k=1,2$ and $d=0$ when the pattern of equal sizes is prescribed, that is, for
fixed $\ast_1,\ast_2,\dots,\ast_{n-1}\in\{=,<\}$, only the sequences satisfying
\[
a_1 \ast_1 a_2 \ast_2 \cdots \ast_{n-1} a_n
\]
are considered.

For a non-decreasing sequence $a=(a_1, a_2, \dots, a_n)$ and $d\ge 0$, let $m(a)$ be the number of rooms with fewer than $a_n-d$ students. We normally shorten $m(a)$ to $m$ if the sequence $a=(a_1,a_2,\dots, a_n)$ is clear from the context. Note that if $d \ge a_n$, then $a$ is
$(k,d)$-dissolvable. Indeed, any two remaining rooms then differ by at most $a_n-1<d$ students, and
this difference stays at most $d$ if we move the students of the dissolved rooms
one by one, each to a remaining room with the fewest students.

Our first main result provides a necessary and sufficient
condition for $(k,d)$-dissolvable sequences.
\begin{theorem}\label{thm:sat}
Let $n,k\in\mathbb{N}$ with $1\le k \le n-1$ and let $d\in \N\cup\{0\}$. Suppose $a =(a_1, a_2, \ldots, a_n)\in \mathbb{N}^n$ is a non-decreasing sequence representing the numbers of students in $n$ classrooms with $a_n > d$. Write $m = \max\{ i : a_i<a_n-d\}$ for the number of rooms with fewer than $a_n-d$ students. Then the following holds.
\begin{enumerate}[label=(\roman*)]
    \item \label{thm:1.1} For $d\ge 1$, $a$ is $(k,d)$-dissolvable if and only if $a_1+ a_2+\dots +a_m\geq (a_n-d)(m-k)$.

    \item \label{thm:1.2} $a$ is $(k,0)$-dissolvable if and only if $a_1+ a_2+\dots +a_n\geq a_n(n-k)$ and $a_1+a_2+\dots +a_n$ is divisible by $n-k$.
\end{enumerate}
\end{theorem}

Intuitively, the hardest case is to dissolve the $k$ smallest rooms, and then every
remaining room must end up with at least $a_n-d$ students since room $n$ remains.
The condition in (i) says exactly that the students in rooms $1,2,\dots,k$ are
enough to bring rooms $k+1,k+2,\dots,m$ up to $a_n-d$ students.
When $d=0$, all remaining rooms must end up with the same number of students, which
is why (ii) also requires $a_1+a_2+\dots+a_n$ to be divisible by $n-k$.

Our second main result provides the values of $N(n,k,d)$ and $N'(n,k,d)$, which have several forms depending on $n,k,$ and $d$.
\begin{theorem}\label{thm:alpha}
Let $n,k\in\mathbb{N}$ with $1\le k \le n-1$ and let $d\in \N\cup\{0\}$. Then we obtain the following:
\begin{enumerate}[label=(\roman*)]
    \item\label{thm2.1} \[
    N(n,k,d) = \begin{cases}
    1 & \text{if} \quad d\ge1, \\
    1 & \text{if} \quad d=0 \text{ and } (n-k)\mid n, \\
    2 & \text{if} \quad d=0 \text{ and }(n-k)\nmid n .
\end{cases}
    \]

    \item\label{thm2.2}
    \[
    N'(n,k,d)= \begin{cases}
    d+\left\lceil \frac{(n-d-1)(n-d)}{2k} \right\rceil& \text{if} \quad 1\le k< \frac{n-d-1}{2}, \\
    n & \text{if} \quad k\geq \frac{n-d-1}{2} \text{ and } d\ge 1, \\
    n & \text{if} \quad k\ge \frac{n-1}{2}, \; d=0 \text{ and } (n-k) \mid \frac{n(n+1)}{2}, \\
    n+1 & \text{if} \quad k\ge \frac{n-1}{2}, \; d=0  \text{ and } (n-k) \nmid \frac{n(n+1)}{2}.
    \end{cases}
    \]
\end{enumerate}

\end{theorem}
Note that $N'(n,k,d)$ grows quadratically in $n$ for fixed $k$ and $d$, while it is
only $n$ or $n+1$ when $k\ge(n-d-1)/2$.
Note also that $N'(n,k,d)=d+N'(n-d,k,0)$ when $k<(n-d-1)/2$, and so allowing a
difference of $d$ has the same effect as removing $d$ classrooms and then adding $d$
students to the largest one.
For $d=0$, \Cref{thm:alpha}\ref{thm2.2} recovers the result of
Sripatak~\cite{Piyashat} stated above.
Since all previous results on this problem, due to Siriseriwan~\cite{Wacharasak},
Sripatak~\cite{Piyashat} and Thamrongthanyalak~\cite{Athipat}, concern the case
$d=0$, \Cref{thm:sat,thm:alpha} are new for $d\ge1$.

The rest of this paper is organized as follows. Section~\ref{sec:thm:sat} is devoted to the proof of Theorem~\ref{thm:sat}. In Section~\ref{sec:thm:alpha}, we prove Theorem~\ref{thm:alpha}.
Finally, we conclude the paper in Section~\ref{sec:conclude} with some open problems.

\section{Proof of Theorem~\ref{thm:sat}}\label{sec:thm:sat}

To prove that $a$ is $(k,d)$-dissolvable, instead of checking every collection of $k$ rooms, it suffices to consider the most difficult scenario, namely the collection consisting of the $k$ smallest rooms. Indeed, dissolving a collection of rooms becomes easier when one of them is replaced by a larger room. The following lemmas make these observations precise.

\begin{lemma}\label{lem:reducecase}
    For $1\le i<i' \le n$, let $I$ be a collection of $k$ rooms including room $i$ but excluding room $i'$ and let $I'$ be obtained from $I$ by replacing room $i$ with room $i'$. If $I$ is $d$-dissolvable, so is $I'$.
\end{lemma}

\begin{proof}
  Since $a$ is non-decreasing, the condition $i<i'$ implies that $a_i\le a_{i'}$.
    To prove that $I'$ is $d$-dissolvable, we first partially dissolve room $i'$ by moving $a_{i'} - a_{i}$ students to room $i$. Then the number of students in room $i$ is equal to $a_{i'}$ and the number of students in room $i'$ is equal to $a_{i}$.
    At this point, the numbers of students are as at the beginning, except that
rooms $i$ and $i'$ have exchanged them.
Therefore, we can finish by applying the dissolution of the rooms in $I$ with
the roles of rooms $i$ and $i'$ exchanged.
    Thus, $I'$ is $d$-dissolvable.
\end{proof}

\begin{lemma}\label{lem:worstcase}
    If rooms $1,2,\dots, k$ are $d$-dissolvable, then every $k$ rooms are $d$-dissolvable.
\end{lemma}

\begin{proof}
    Let $1\le i_1<i_2<\dots< i_k \le n$. Note that $i_j\geq j$ for all $j$.
    Since rooms $1,2,\dots, k$ are $d$-dissolvable, by applying \Cref{lem:reducecase}, so are rooms $1,2,\dots, k-1, i_k$. Next, by applying \Cref{lem:reducecase} again, rooms $1,2,\dots,k-2,i_{k-1}, i_k$ are $d$-dissolvable. Proceeding inductively, we
obtain that rooms $i_1,i_2,\dots, i_k$ are $d$-dissolvable.
\end{proof}
The following lemma gives the necessary condition in \Cref{thm:sat}
\ref{thm:1.1}, and we state it for every $d\ge0$ since we shall also use it when
$d=0$.
\begin{lemma}\label{lem:morethand}
Let
$d\ge0$. If $a_n >d$ and $a$ is $(k, d)$-dissolvable, then $a_1 +\dots+ a_m\geq (a_n-d)(m- k)$.
\end{lemma}
\begin{proof}
   Note that the inequality always holds if $m<k$. Throughout this proof, we suppose that $m\ge k$. Assume that $a$ is $(k,d)$-dissolvable. Then rooms $1,2,\dots, k$ are $d$-dissolvable. After dissolving, since the difference between the numbers of students in any two rooms is at most $d$ and room $n$ is not dissolved, the size of each remaining room is at least $a_n-d$.
    For room $i = k+1,k+2,\dots, m$, we must have moved at least $a_n-d-a_i$ students to room $i$.
    Hence, the total number of students from the dissolved rooms $a_1+a_2+\dots+a_k$ must be at least $\sum_{i=k+1}^{m} (a_n-d-a_{i})$, i.e.\ $a_1+ a_2+\dots +a_m\geq (a_n-d)(m-k)$.
\end{proof}
We are now ready to prove \Cref{thm:sat}.

\begin{proof}[Proof of \Cref{thm:sat}\ref{thm:1.1}]
First, the forward direction follows from \Cref{lem:morethand}.

Conversely, assume that $a_1+ a_2+\dots +a_m\geq (a_n-d)(m-k)$.
    To show that $a$ is $(k,d)$-dissolvable, by using \Cref{lem:worstcase}, it suffices to show that the smallest rooms $1,2,\dots, k$ are $d$-dissolvable. Since $a_1+a_2+\dots+a_k \ge \sum_{i=k+1}^{m} (a_n-d-a_{i})$, which means that the total number of students in the dissolved rooms is sufficient to compensate for those rooms whose sizes are less than $a_n-d$, we move $\sum_{i=k+1}^{m} (a_n-d-a_{i})$ students to rooms $k+1, k+2, \dots, m$, raising their sizes to $a_n-d$.
     After that, we distribute the remaining $(a_1+a_2+\dots+a_k)-\sum_{i=k+1}^{m}
(a_n-d-a_i)$ students among the $n-k$ remaining rooms so that the numbers of
students in these rooms are as equal as possible. As a result, since $d \ge 1$, the difference between the numbers of students in any two rooms is at most $d$. Hence, the sequence $a$ is $(k,d)$-dissolvable.
\end{proof}

\begin{proof}[Proof of \Cref{thm:sat}\ref{thm:1.2}]
  First, we get from \Cref{lem:morethand} that $a_1+ a_2+\dots +a_m\geq a_n(m-k)$.
    Adding $a_n(n- m) $ to both
sides, we obtain $a_1+ a_2+\dots +a_n\geq a_n(n-k)$. Since the remaining $n-k$ rooms after rearrangement have the same size, the total number of students, which is equal to $a_1+a_2+\dots +a_n$, is divisible by $n-k$.

    Conversely, assume that $a_1+ a_2+\dots +a_n\geq a_n(n-k)$ and $a_1+a_2+\dots+a_n$ is divisible by $n-k$.
    To show that $a$ is $(k,0)$-dissolvable, by using \Cref{lem:worstcase}, it suffices to show that the smallest rooms $1,2,\dots, k$ are $0$-dissolvable. Since $a_1+a_2+\dots+a_k \ge \sum_{i=k+1}^{n} (a_n-a_{i})$, we move $\sum_{i=k+1}^{n} (a_n-a_{i})$ students to rooms $k+1, k+2, \dots, n-1$, raising their sizes to $a_n$.
    After that, we distribute the remaining $(a_1+a_2+\dots+a_k) - \sum_{i=k+1}^{n} (a_n-a_{i})$ students as equally as possible among the $n-k$ remaining rooms.
    Since the total number of students is divisible by $n-k$, all remaining rooms have the same size. Hence, the sequence $a$ is $(k,0)$-dissolvable.
\end{proof}

\section{Proof of Theorem~\ref{thm:alpha}}\label{sec:thm:alpha}
To determine the values of $N(n,k,d)$ and $N'(n,k,d)$, we first derive a lower bound for each quantity and then construct examples showing that the bounds are sharp. Applying \Cref{thm:sat}, we verify
that the examples are dissolvable. We begin by determining $N(n,k,d)$.
\begin{proof}[Proof of \Cref{thm:alpha}\ref{thm2.1}]
    If $d\ge1$, then the sequence $(1,1,\dots,1)$ has $a_n=1\le d$, and so it is
$(k,d)$-dissolvable as observed in Section 1. Hence, $N(n,k,d)=1$. Suppose now that $d=0$. For $(n-k)\mid n$, the sequence $(1,1,\dots, 1)$ has $n$
students in total, which is greater than $n-k$, so the sequence $(1,1,\dots, 1)$ is $(k,0)$-dissolvable by \Cref{thm:sat}\ref{thm:1.2}. Thus, $N(n,k,0) =1$ for $(n-k)\mid n$.
    For $(n-k)\nmid n$,
the sequence $(1,1,\dots, 1)$ is not $(k,0)$-dissolvable by \Cref{thm:sat}\ref{thm:1.2}. Thus, $N(n,k,0) \geq2$ for $(n-k)\nmid n$.%

Next, we will show that $N(n,k,0) =2$ provided $(n-k) \nmid n$. We consider two cases, $2k<n$ and $n<2k$.
    For the case $2k<n$, the sequence consisting of $2k$ ones and $n-2k$ twos has $2(n-k)=a_n(n-k)$
students in total, and so it is $(k,0)$-dissolvable by \Cref{thm:sat}\ref{thm:1.2}. Now, we prove the case $n<2k$. Choose $r$ as the remainder when $2n$ is divided by $n-k$. Then we consider the sequence consisting of $r$ ones and $n-r$ twos. This sequence has $2n-r$ students in total, which is divisible by $n-k$ and
greater than $a_n(n-k)=2(n-k)$ since $r<n-k<k$. Hence, it is $(k,0)$-dissolvable by \Cref{thm:sat}\ref{thm:1.2}.
\end{proof}

We now turn to the determination of $N'(n,k,d)$.

\begin{proof}[Proof of \Cref{thm:alpha}\ref{thm2.2}] Let $d$ be a non-negative integer and $a = (a_1,a_2,\dots, a_n)$ a strictly increasing $(k,d)$-dissolvable sequence.
Then $a_n>a_{n-1}>\dots>a_1\ge 1$, so $a_n\ge n$, which implies $N'(n,k,d)\ge n$. We divide the rest of the proof into four cases, according to the four lines of
the formula in \Cref{thm:alpha}\ref{thm2.2}.

\noindent \textbf{Case 1:} $1\le k< (n-d-1)/2$. We first derive the lower bound of $N'(n,k,d)$.
Since $a_n\ge n> d+2k+1 >d$, we obtain from \Cref{lem:morethand} that
\begin{equation}\label{eq:thm1cond}
    a_1+a_2+\cdots +a_{m} \ge (a_n-d)(m-k).
\end{equation}
Since $a$ is strictly increasing, i.e., $a_1<a_2<\dots <a_n$, it follows that $a_i \le a_n-(n-i)$ for any $i = 1,2,\dots, n$.
Then $a_1<a_2<\dots<a_{n-d-1} \le a_n-d-1 < a_n-d$. From $m = \max\{i : a_i<a_n-d\}$, we have
\begin{equation}\label{eq:mbound}
    n-d-1 \le m \le n-1.
\end{equation}
Since $a_1<a_2<\dots<a_m<a_n-d$, we obtain
\begin{equation}\label{eq:aibound}
    a_i \le a_n-d-(m+1-i) \quad  \text{for each } i=1, 2, \ldots, m.
\end{equation}
Then
\[
(a_n-d)(m-k) \overset{\eqref{eq:thm1cond}}{\le} \sum_{i=1}^{m} a_i \overset{\eqref{eq:aibound}}{\le} \sum_{i=1}^{m} (a_n-d-m-1+i) = m(a_n-d) - \frac{m(m+1)}{2},
\]
and so,
\[
a_n \ge d + \frac{m(m+1)}{2k} \overset{\eqref{eq:mbound}}{\ge} d + \frac{(n-d-1)(n-d)}{2k}.
\]
Since $a_n \in \N$, we have
\[
a_n \ge d + \left\lceil \frac{(n-d-1)(n-d)}{2k} \right\rceil.
\]
Thus, $N'(n,k,d) \ge d + \left\lceil \frac{(n-d-1)(n-d)}{2k} \right\rceil$.

Next, we construct examples attaining the lower bound for each $d$. Since \eqref{eq:mbound} and \eqref{eq:aibound} hold with equality when the
sequence consists of consecutive integers, our examples will be consecutive
integers ending at the lower bound, adjusted slightly when $d=0$. For $d\ge 1$, consider a sequence $b = (b_1, b_2, \ldots, b_n)$ with $b_n = d + \left\lceil \frac{(n-d-1)(n-d)}{2k} \right\rceil$ and $b_i = b_n - (n - i)$ for $1 \le i \le n-1$.
It follows from $2k<n-d-1$ that $b_n > n$, and so, $b_1 = b_n-(n-1)>1$. Note that $m(b)=n-d-1$, since $b_i<b_n-d$ exactly when $i\le n-d-1$.
We also have
\begin{align*}
    b_1 + b_2 + \ldots + b_{m} &= mb_n - mn + \frac{m(m+1)}{2} \\
    &= m(b_n-d) + m \bigg( \frac{m+1}{2} - (n-d) \bigg) \\
    &= m(b_n-d) + (n-d-1) \bigg( \frac{n-d}{2} - (n-d) \bigg) \\
    &= (m-k)(b_n-d) + k \left\lceil \frac{(n-d-1)(n-d)}{2k} \right\rceil - \frac{(n-d-1)(n-d)}{2} \\
    &\ge (m-k)(b_n-d).
\end{align*}
Hence, $b$ is $(k,d)$-dissolvable by \Cref{thm:sat}\ref{thm:1.1}, and thus, $N'(n,k,d) = d + \left\lceil \frac{(n-d-1)(n-d)}{2k} \right\rceil$.

For $d=0$, observe that the sum of $n$ consecutive integers ending at $\left\lceil \frac{(n-1)n}{2k} \right\rceil$ is
\begin{equation}\label{eq:N'd=0}
    \sum_{i=1}^n \left( \left\lceil \frac{(n-1)n}{2k} \right\rceil-(n-i) \right) =(n-k)\left\lceil \frac{(n-1)n}{2k} \right\rceil + k\left\lceil \frac{(n-1)n}{2k} \right\rceil - \frac{(n-1)n}{2}.
\end{equation}
To make the sum divisible by $n-k$, let $r = k\left\lceil \frac{(n-1)n}{2k} \right\rceil - \frac{(n-1)n}{2} \in\Z$. Note that $r\ge 0$. Consider $c = (c_1,c_2,\dots, c_n)$ with $c_n = \left\lceil \frac{(n-1)n}{2k} \right\rceil$ and
\[
c_i = \begin{cases}
    c_n-(n-i)-1 & \text{if} \quad 1\le i \le r, \\
    c_n-(n-i) & \text{if} \quad r+1 \le i \le n.
\end{cases}
\]
It follows from $2k<n-1$ that $c_n > n$, and so, $c_1 \geq c_n-n>0$, which implies $c_1\ge 1$. Then
\[
\sum_{i=1}^n c_i = \sum_{i=1}^n \big( c_n-(n-i) \big) - r \overset{\eqref{eq:N'd=0}}{=} c_n(n-k).
\]
By \Cref{thm:sat}, $c$ is $(k,0)$-dissolvable, and hence, $N'(n,k,0) = \left\lceil \frac{(n-1)n}{2k} \right\rceil$.

\noindent\textbf{Case 2:} $2k\ge n-d-1$ and $d\ge 1$.
Consider $a=(a_1, a_2, \ldots, a_n)=(1,2,\dots, n)$. If $d\geq n$, then $a$ is $(k,d)$-dissolvable as observed in Section 1. If $d<n$, we have $m = n-d- 1 \geq 0$ and $m\le 2k$. Thus,
\[
a_1+a_2+\dots+a_m = \frac{m(m+1)}{2} = (m+1)\left( m-\frac{m}{2} \right) \ge (n-d)(m-k).
\]
It follows from \Cref{thm:sat} that $a$ is $(k,d)$-dissolvable, and thus, $N'(n,k,d)=n$.

\noindent\textbf{Case 3:} $2k\ge n-1, \; d=0$ and $(n-k) \mid \frac{n(n+1)}{2}$.
Consider $a=(a_1, a_2, \ldots, a_n)=(1,2,\dots, n)$. Since $2k\ge n-1$, we have $(n+1)/2\geq n-k$. Then $a_1+a_2+\dots+a_n=n(n+1)/2\geq n(n-k)=a_n(n-k)$ and $(n-k)\mid \frac{n(n+1)}{2}.$ It follows from \Cref{thm:sat} that $a$ is $(k,0)$-dissolvable, and thus, $N'(n,k,0)=n$.

\noindent \textbf{Case 4:} $2k\ge n-1, \; d=0$ and $(n-k) \nmid \frac{n(n+1)}{2}$.
It is clear that the sequence $(1,2,\dots, n)$ is not $(k,0)$-dissolvable, so $N'(n,k,0)\ge n+1$.
 Let $r\in\{1,2,\dots,n-k\}$ with $r\equiv (n+1)(n+2)/2 \pmod{n-k}$, and let
$a$ be the sequence obtained from $(1,2,\dots,n+1)$ by removing $r$.
We have $a_1+a_2+\dots+ a_n=1+2+\dots+n+(n+1)-r = (n+1)(n+2)/2-r$. Then $(n-k)\mid a_1+a_2+\dots+a_n$ and
\[
a_1+a_2+\dots+a_n=(n+1)(n+2)/2-r\ge (n+1)(n+2)/2-(n-k)\geq (n+1)(n-k).
\]
It follows from \Cref{thm:sat} that $a$ is $(k,0)$-dissolvable, and thus, $N'(n,k,0)=n+1.$
\end{proof}

\section{Concluding remarks}
\label{sec:conclude}
We establish a necessary and sufficient condition for $(k,d)$-dissolvability and determine the exact values of $N(n,k,d)$ and $N'(n,k,d)$ for all $n,k$ and $d$.
Note that $N(n,k,d)$ allows all classroom sizes before redistribution to be equal, while $N'(n,k,d)$ requires them to be pairwise distinct. By \Cref{thm:sat} and the observation before it, the $(k,d)$-dissolvable sequences
are characterized for all $k$ and $d$, including the sequences considered
in~\cite{Athipat}. However, for fixed $\ast_1,\ast_2,\dots,\ast_{n-1}\in\{=,<\}$,
the minimum value of $a_n$ over all $(k,d)$-dissolvable sequences satisfying
$a_1\ast_1 a_2\ast_2\cdots\ast_{n-1}a_n$ is not known in general for $k\ge3$ or
$d\ge1$.
Furthermore, we consider rearrangements for which the numbers of students in the classrooms differ pairwise by at most $d$, with no further restrictions. It would be interesting to investigate analogous problems in which the classroom sizes after regrouping are required to satisfy additional conditions, for example, forming an arithmetic progression or requiring some classrooms to have equal sizes.

\section*{Acknowledgements}
The research of this paper was carried out while the authors were participating in the Geometry and Combinatorics Boot Camp 2025, in Bangkok, Thailand.
We are grateful for the hospitality of the Department of Mathematics and Computer Science, Faculty of Science, Chulalongkorn University.

\bibliographystyle{siam}
\bibliography{partsat}

\end{document}